\documentclass[journal,twoside,web]{ieeecolor}

\usepackage{generic}
\usepackage{cite}
\usepackage{amsmath,amsthm,amssymb,amsfonts}
\usepackage{algorithmic}
\usepackage{graphicx}
\usepackage{textcomp}
\usepackage{subfigure}
\usepackage{epsfig}
\usepackage{lcsys}

\newtheorem{theorem}{Theorem}[section]

\newtheorem{lemma}[theorem]{Lemma}
\newtheorem{proposition}[theorem]{Proposition}
\newtheorem{corollary}[theorem]{Corollary}
\newtheorem{definition}[theorem]{Definition}
\newtheorem{remark}[theorem]{Remark}
\newtheorem{assumption}[theorem]{Assumption}

\usepackage{amsmath,amssymb,amsthm,mathtools,bm}
\usepackage{microtype}
\let\labelindent\relax 
\usepackage{enumitem}
\usepackage{booktabs}
\usepackage[T1]{fontenc}
\usepackage{lmodern}

\newcommand{\R}{\mathbb{R}}
\newcommand{\cI}{\mathrm{ci}}
\newcommand{\rci}{\mathrm{rci}}
\newcommand{\st}{\mathrm{s.t.}}
\newcommand{\conv}{\operatorname{conv}}
\newcommand{\vtx}{\operatorname{vert}}
\newcommand{\one}{\bm 1}
\newcommand{\rev}[1]{\textcolor{black}{#1}}   
\newcommand{\mrev}[1]{{\color{black}#1}}     
\newcommand{\revtwo}[1]{\textcolor{black}{#1}}  
\definecolor{revpurple}{RGB}{128,0,128}
\newcommand{\revthree}[1]{\textcolor{black}{#1}} 

\newcounter{algorithm}

\begin{document}

\def\BibTeX{{\rm B\kern-.05em{\sc i\kern-.025em b}\kern-.08em
    T\kern-.1667em\lower.7ex\hbox{E}\kern-.125emX}}
\markboth{IEEE Control Systems Letters, VOL. XX, NO. XX, 2026}
{Athanasopoulos \MakeLowercase{\textit{et al.}}: Safety Certificates via Stackelberg Set Games}

\title{Safety Certificates for Systems under Attack \\ via Stackelberg Set Games}

\author{Nikolaos Athanasopoulos, David Laverty, Gregory Provan, and Quanyan Zhu
\thanks{The authors are listed in alphabetical order. N. A. is with the Department of Electrical and Computer Engineering, University of Thessaly, Greece (e-mail: \texttt{n.athanasopoulos@uth.gr}), D. L. is with the School of Electronics, Electrical Engineering and Computer Science, Queen's University Belfast, Belfast, U.K. (e-mail:  \texttt{david.laverty@qub.ac.uk}). G. P. is with the School of Computer Science and Information Technology, University College Cork, Cork, Ireland (e-mail: \texttt{g.provan@ucc.ie}). Q. Z. is with the Department of Electrical and Computer Engineering, New York University, Brooklyn, NY 11201 USA (e-mail: \texttt{qz494@nyu.edu}). All authors gratefully acknowledge support from the U.S.-Ireland R\&D Partnership under NSF Grant \rev{No.~2521130}.}}

\IEEEaftertitletext{\vspace{-3\baselineskip}}
\maketitle
\pagestyle{empty}
\thispagestyle{empty}

\begin{abstract}
We study a Stackelberg set selection game for linear systems subject to
bounded disturbances and additive attacks. The defender selects the
vector parameter of a polytopic safe set template; the attacker, after
observing this choice, selects the vector parameter of an admissible
attack set. The payoff is a signed, facet-normalised invariance margin
that tests invariance exactly and quantifies either the available
contraction or the facet inflation needed when invariance fails. Both controlled and robust controlled invariance are treated,
corresponding to different information patterns for the control input.
Exploiting the structure of parameterised polyhedra, we show that within each validity
region the margin is concave in the attack parameter, which reduces the
follower minimisation to finitely many extreme parameters. The complete max--min game is then solved by a finite family of LPs, one
per defender validity region, each enforcing all attacker regions under a
common defender decision and guaranteed margin. The framework handles
hard state constraints and additional costs related to the safe set size, and
recovers optimal set parameters together with state and attack
realisations that activate the certificate. A load frequency control example under false data
injection illustrates the construction, the validity decomposition, and
the two information patterns.
\end{abstract}


\section{Introduction}
\label{sec:novelty}

\revtwo{
\revthree{Safety of a constrained control system can be certified by
constructing an invariant set, namely,} a region of the state space within which admissible control inputs keep the system trajectory despite disturbances \cite{Blanchini1999}. Motivated by adversarial control settings \cite{Cardenas2008,Fawzi2014,Teixeira2015,ZhuBasar2015}, we consider an attacker who observes the defensive design and selects the set of admissible additive attacks subject to prescribed budget constraints. The resulting design question is which safe set the defender should
select when anticipating the attacker’s response.
}

Set invariance and polyhedral vertex based certificates are classical tools
in constrained control \cite{Blanchini1999}. The optimisation of
invariant polytopes with prescribed complexity or facet directions is
also established \cite{Athanasopoulos2014,Trodden2016}, including joint
facet and vertex parameterisations \cite{Villanueva2024}. The affine
vertex decompositions used in this note follow from parameterised polyhedra
 \cite{LoechnerWilde1997}; related decompositions have also been
used for state dependent polytopic disturbances
\cite{Athanasopoulos2025}. The leader--follower order and its
information pattern are standard in Stackelberg games
\cite{BasarOlsder1999}\rev{, and the resulting max--min programmes are
bilevel in general \cite{SinhaMaloDeb2018}}.

\revtwo{While none of the above individual
ingredients is new, our contribution is their coupling in a \emph{Stackelberg game over sets}}:
The players choose the shape of polytopes within predefined templates, rather than a single state, an
attack realisation, or a scalar uncertainty dilation. The payoff is
a signed, facet normalised invariance margin for two control information
patterns. Despite continuous set parameters and
piecewise affine vertex maps,  \rev{the follower parameter is eliminated
exactly, while a single safe set offset and a single guaranteed margin
are retained across all attacker regions}. The complete game is therefore solved by one LP per
defender validity region. 
Thus, our key contribution is the observation that strategic set selection
preserves the LP structure for polytopic invariance,
with complexity governed by the validity regions and generators.

Three 
observations deserve emphasis: First, the signed margin turns
a binary invariance test into an objective: the sign decides invariance
exactly, while magnitude measures available contraction and
required inflation on a single facet normalised scale, thus safety
can be traded against safe set size in a design problem. \revtwo{Second, for each defender parameter, the margin is concave in the attack parameter within a validity region. Consequently, an attacker best response can be selected from the union of the vertices of these regions. This finite reduction is exact, not requiring discretising the admissible parameter set}.
Third, the order of play cannot be recovered from the individual
regional problems, because the defender commits to one safe set before
the attacker selects its validity region, so optimising those problems
separately solves a different and generally more optimistic problem.
Compared with robust and reachability based certificates, where the
adversary perturbs a realisation inside a fixed uncertainty set, the
adversary here selects the uncertainty set itself subject to a
budget.

\section{Preliminaries and setting}
\label{sec:polytopic-game}

\emph{Notation:} 
\revthree{We consider a norm $\|\cdot\|$ with unit ball $\mathbb B$ and
its dual $\|\cdot\|_*$.} Vector inequalities hold componentwise. \rev{The vector of all ones is $\one$, and $\vtx(\cdot)$ denotes the
vertex set of a polytope.}

The focus \rev{of} this note is on linear dynamics and polytopic sets. Starting from a state constraint set, the defender selects the shape of a candidate safe set template, while the attacker selects the
offsets of an admissible attack set. 
We consider the attacked system
\begin{equation}
 x^+=Ax+Bu+c+Ed+\delta,
 \label{eq:system}
\end{equation}
where $x\in\R^n$, $u\in U\subseteq\R^{n_u}$ is the defender control \rev{input}, $d\in D$ is the nondeterministic additive disturbance, and
$\delta\in\R^n$ is the
additive attack vector. The control and disturbance sets are the compact polytopes
\begin{equation}
 U=\{u:F_u u\leq g_u\},\qquad
 D=\conv\{d_1,\ldots,d_{N_D}\}.
 \label{eq:UD}
\end{equation}
\revthree{Here $F_u$ has $m_U$ rows; state constraints are}
$x\in\mathcal{X}$, with  
\begin{equation}
 \mathcal{X}:=\{x:H_Xx\leq h_X\}.
 \label{eq:hard-X}
\end{equation}
\revthree{The defender and attacker sets are $S(b)$ and $\Delta(\rho)$,}
with
\begin{align}
 S(b)&:=\{x\in\R^n:Cx\leq b\},
 &&b\in\mathcal B\subseteq\R^p, \label{eq:Sb}\\
 \Delta(\rho)&:=\{\delta\in\R^n:C_\Delta\delta\leq\rho\},
 &&\rho\in\mathcal R\subseteq\R^{p_\Delta},
 \label{eq:Deltarho}
\end{align}
The strategic variables $b,\rho$ parameterise sets instead of individual
states or attack realisations. We note the admissible set $\mathcal B$ can
incorporate $S(b)\subseteq \mrev{\mathcal X}$ directly, or this
specification \rev{can be} handled explicitly in a design problem, as
\rev{outlined in Section~\ref{sec:design}}.
\begin{assumption}\label{ass:polytopes}
\revthree{The sets $\mathcal B$, $\mathcal R$, $U$ and $D$ are nonempty
compact polytopes} for each admissible pair $(b,\rho)$,
and the same holds for $S(b)$ and $\Delta(\rho)$.
\end{assumption}

The families in \eqref{eq:Sb} and \eqref{eq:Deltarho} are special cases of 
parameterised polyhedral sets, since their facet normals are fixed while the offset varies affinely on the parameters. The finite parameterised vertex representation in
\cite{LoechnerWilde1997}, explicitly stated in \cite{Athanasopoulos2025}, yields a finite polytopic partition 
\begin{equation}
 \mathcal B=\bigcup_{j\in\mathcal J}\mathcal B_j,
 \qquad
 \mathcal R=\bigcup_{h\in\mathcal H}\mathcal R_h.
 \label{eq:partitions}
\end{equation}
Specifically, for each $j\in \mathcal{J}$, $h\in\mathcal{H}$ we have
\begin{align}
 S(b)&=\conv\{v_1^j(b),...,v_{q_j}^j(b)\},
 \ v_r^j(b)=M_r^j b+m_r^j, \label{eq:Svertices}\\
 \Delta(\rho)&=\conv\{w_1^h(\rho),...,
 w_{q_{\Delta,h}}^h(\rho)\}, 
 w_s^h(\rho)=N_s^h\rho+n_s^h. \label{eq:Dvertices}
\end{align}
Each defender region admits the representation
$\mathcal B_j=\{b:G_jb\leq\beta_j\}$ with
$G_j\in\R^{m_{B,j}\times p}$. We discard empty validity domains and \rev{allow different domains to
overlap on their boundaries. The vertex representation is piecewise
affine globally, however it is affine on any single region}
$\mathcal{B}_j$, $\mathcal{R}_h$.
We are interested in two notions of controlled invariance.
\begin{definition} 
\label{def:invariance}
Consider \eqref{eq:system} with $U$, $D$ as in \eqref{eq:UD} and
$S(b)$, $\Delta(\rho)$ as in
\eqref{eq:Sb}--\eqref{eq:Deltarho}. For fixed
$(b,\rho)\in\mathcal B\times\mathcal R$, the set $S(b)$ is:
(i) \emph{controlled invariant} if $\forall x\in S(b)$,
$\forall d\in D$, $\forall\delta\in\Delta(\rho)$,
$\exists u\in U:x^+\in S(b)$; (ii) \emph{robust controlled invariant}
if $\forall x\in S(b)$, $\exists u\in U$, $\forall d\in D$,
$\forall\delta\in\Delta(\rho):x^+\in S(b)$.
\end{definition}
Consequently, for controlled invariance,  the input depends on $(x,d,\delta)$, while
for robust controlled invariance  the input depends on $x$ alone. In both cases, the selected parameters $(b,\rho)$ are
known to the controller; the distinction concerns access
to the realised disturbance and attack. \revthree{Without changing $S(b)$, we normalise the nonzero rows of $C$
so that $\|C_i^\top\|_*=1$, rescaling the corresponding components of
$b$.} For $\varepsilon\in\R$ we define the safety margin set
$S_\varepsilon(b):=\{z:Cz\leq b-\varepsilon\mrev{\one}\}$.
For $\varepsilon\geq0$, $S_\varepsilon(b)=S(b)\ominus\varepsilon\mathbb B$;
for $\varepsilon<0$ it is a facetwise outward relaxation. In general,
$S(b)\oplus r\mathbb B\subseteq S_{-r}(b)$ for $r=-\varepsilon>0$,
thus a negative margin should not be read
as a Minkowski enlargement.

\begin{definition}
In the setting of Definition~\ref{def:invariance}, 
let
$\nu\in\{\cI,\rci\}$. The \emph{signed invariance margins}
$\varepsilon_\nu:\mathcal B\times\mathcal R\to\R$ are
\begin{align*}
 \varepsilon_{\cI}(b,\rho):=\max\{\varepsilon\in\R:\;&
 \forall x\in S(b),\ \forall d\in D,
 \ \forall\delta\in\Delta(\rho),\nonumber\\[-1mm]
 &\exists u\in U:\ x^+\in S_\varepsilon(b)\},
 \\
 \varepsilon_{\rci}(b,\rho):=\max\{\varepsilon\in\R:\;&
 \forall x\in S(b),\ \exists u\in U,
 \ \forall d\in D,\nonumber\\[-1mm]
 &\forall\delta\in\Delta(\rho):\ x^+\in S_\varepsilon(b)\}.
\end{align*}
\end{definition}
For fixed $(b,\rho)$ the feasible margin set is nonempty and bounded
above, since compactness bounds $x^+$; closedness of $S(b)$ and
compactness of $U$ make both maxima attained.
The sign of $\varepsilon_\nu$ is an
exact invariance test:
$\varepsilon_\nu(b,\rho)\geq0$ if and only if $S(b)$ satisfies $\nu$.
A negative margin implies invariance fails and quantifies the required
relaxation; a nonnegative one certifies invariance and measures the
available contraction. Being scalar, $\varepsilon_\nu$ supplies the total ordering a zero sum payoff needs.
We note the robust pattern admits fewer control choices, i.e.,
$\varepsilon_{\rci}(b,\rho)\leq\varepsilon_{\cI}(b,\rho)$ for any pair
$(b,\rho)\in\mathcal B\times\mathcal R$. Thus, any nonnegative robust
controlled invariance margin also certifies controlled invariance,
whereas the converse does not hold in general.

\section{Main result}
\revthree{We derive finite programmes for the invariance margins and the
associated Stackelberg set game: exact vertex certificates at fixed set
parameters, a reduction of the attacker minimisation to finitely many
extreme parameters one LP per defender region enforcing all attacker regions under a common
parameter and margin, and the global solution by comparing their
values.}
\begin{proposition}[Exact regional certificates]
\label{prop:regional-certificates}
Let $j\in\mathcal J$ and $h\in\mathcal H$. For any pair $(b,\rho)\in\mathcal{B}_j\times \mathcal{R}_h$, the controlled
invariance margin $\varepsilon_{\mrev{\cI}}(b,\rho)$ is the optimal value of
\begin{subequations}\label{eq:ci-fixed}
\begin{align}
 \max_{\varepsilon,\{u_{r,l,s}\}}\quad&\varepsilon
 \\
 \st\quad
 &C\!\left(Av_r^j(b)+Bu_{r,l,s}+c\right)
 \nonumber\\[-1mm]
 &\quad+C\!\left(Ed_l+w_s^h(\rho)\right)
 \leq b-\varepsilon\mrev{\one},
 \quad\forall r,l,s,
 \\
 &F_u u_{r,l,s}\leq g_u,
 \hspace{14mm}\forall r,l,s.
\end{align}
\end{subequations}
Respectively, the robust controlled invariance margin $\varepsilon_{\text{rci}}(b,\rho)$ is the optimal value of
 \begin{subequations}\label{eq:rci-fixed}
\begin{align}
 \max_{\varepsilon,\{u_r\}}\quad&\varepsilon
 \\
 \st\quad
 &C\!\left(Av_r^j(b)+Bu_r+c\right)
 \nonumber\\[-1mm]
 &\quad+C\!\left(Ed_l+w_s^h(\rho)\right)
 \leq b-\varepsilon\mrev{\one},
 \quad\forall r,l,s,
 \\
 &F_u u_r\leq g_u,
 \hspace{22mm}\forall r.
\end{align}
\end{subequations}
In \eqref{eq:ci-fixed}--\eqref{eq:rci-fixed}, $r=1,\ldots,q_j$,
$l=1,\ldots,N_D$, $s=1,\ldots,q_{\Delta,h}$.
\end{proposition}

\begin{proof}
To show necessity, we apply the corresponding quantified condition to
each vertex triple $(v_r^j(b),d_l,w_s^h(\rho))$, precisely leading to \eqref{eq:ci-fixed} or \eqref{eq:rci-fixed}.
For sufficiency, we observe that any
$(x,d,\delta)\in S(b)\times D\times\Delta(\rho)$ is a convex
combination of the vertex triples\rev{, with coefficients
$\alpha_r\beta_l\gamma_s$ formed as products of the barycentric
coefficients of $x$, $d$ and $\delta$, whose marginals reproduce the
three components}. For controlled invariance, the same
product coefficients combine the associated inputs $u_{r,l,s}$. For
robust controlled invariance, writing
$x=\sum_r\alpha_rv_r^j(b)$ and selecting
$u=\sum_r\alpha_ru_r$ gives an admissible input that is independent of
the subsequently realised $(d,\delta)$. Convexity of $U$, affinity of
\eqref{eq:system}, and convexity of $S_\varepsilon(b)$ establish
the corresponding quantified condition.\end{proof}

\revthree{In the sequel, after fixing $\mathcal B_j$ and $\mathcal R_h$,
we replace the inequalities quantified by $\rho$ with finitely many
inequalities, for every $b$.} 
\rev{Throughout, $\varepsilon_\nu^{j,h}(b,\rho)$ denotes the optimal value of the \revthree{programme of
Proposition~\ref{prop:regional-certificates} with the affine maps of
$(j,h)$}. It coincides with $\varepsilon_\nu(b,\rho)$ for every
$(b,\rho)\in\mathcal B_j\times\mathcal R_h$, so \revthree{parameters in
overlapping domains receive the same value}. The
\emph{regional game value} is}
\begin{equation}
 V_\nu^{j,h}:=
 \max_{b\in\mathcal B_j}
 \min_{\rho\in\mathcal R_h}
 \varepsilon_\nu^{j,h}(b,\rho),
 \qquad \nu\in\{\cI,\rci\},
 \label{eq:regional-game}
\end{equation}
with the superscript indicating that the affine vertex maps
\eqref{eq:Svertices}--\eqref{eq:Dvertices} are fixed.
Since $\mathcal R_h$ is a compact polytope, its vertex set
$\vtx(\mathcal R_h)=\{\rho^{h,1},\ldots,\rho^{h,K_h}\}$ is finite,
with $K_h=|\vtx(\mathcal R_h)|$.
\revthree{The following result writes the max--min problem
\eqref{eq:regional-game} as a maximisation with finitely many hypograph
inequalities; Proposition~\ref{prop:regional-LP} then replaces each by linear constraints, producing one LP.}
\begin{lemma}[Concavity and finite hypograph representation]
\label{lem:finite-hypograph}
For every $j\in\mathcal J$, $h\in\mathcal H$, and
$b\in\mathcal B_j$, the mapping
$\varepsilon_\nu^{j,h}(b,\cdot)$ is concave on $\mathcal R_h$.
Moreover,
\begin{subequations}\label{eq:regional-hypograph}
\begin{align}
 V_\nu^{j,h}=\max_{b,t}\quad&t
 \\
 \st\quad&b\in\mathcal B_j,
 \\
 &t\leq\varepsilon_\nu^{j,h}(b,\rho^{h,k}),
 \qquad k=1,\ldots,K_h.
 \label{eq:regional-hypograph-con}
\end{align}
\end{subequations}
\end{lemma}

\begin{proof}
Take $t$ as a common lower bound on all attacker responses, so that for
$b\in\mathcal B_j$
\begin{align}
 \min_{\rho\in\mathcal R_h}\varepsilon_\nu^{j,h}(b,\rho)
 =\max_t\{t:\;       t\leq\varepsilon_\nu^{j,h}(b,\rho),
 \forall\rho\in\mathcal R_h\}.
 \label{eq:hypograph-identity}
\end{align}
For concavity, fix $b\in\mathcal B_j$, take
$\rho^1,\rho^2\in\mathcal R_h$ and $\theta\in[0,1]$, and let
$(\varepsilon^a,\bm u^a)$ optimise the relevant LP of
Proposition~\ref{prop:regional-certificates} at $\rho^a$. Set
$\rho^\theta:=\theta\rho^1+(1-\theta)\rho^2$,
$\varepsilon^\theta:=\theta\varepsilon^1+(1-\theta)\varepsilon^2$,
and combine the controls in the same way,
$u_{r,l,s}^\theta$ for controlled invariance and $u_r^\theta$ for the
robust pattern; both lie in $U$ by convexity. Since $w_s^h(\cdot)$ is
affine on $\mathcal R_h$, every $x^+$ constraint at
$(\rho^\theta,\varepsilon^\theta,\bm u^\theta)$ is the convex
combination of those at $(\rho^a,\varepsilon^a,\bm u^a)$ thus it is
satisfied, so
$\varepsilon_\nu^{j,h}(b,\rho^\theta)\geq\varepsilon^\theta$, which
is concavity. Writing $\rho\in\mathcal R_h$ as
$\sum_k\theta_k\rho^{h,k}$ with $\theta\geq0$,
$\sum_k\theta_k=1$, concavity gives
$\varepsilon_\nu^{j,h}(b,\rho)\geq
\sum_k\theta_k\varepsilon_\nu^{j,h}(b,\rho^{h,k})
\geq\min_k\varepsilon_\nu^{j,h}(b,\rho^{h,k})$, so the minimum over
$\mathcal R_h$ is at a vertex. Substituting in
\eqref{eq:hypograph-identity} 
gives \eqref{eq:regional-hypograph}.
\end{proof}

\revthree{The $K_h$ inequalities \eqref{eq:regional-hypograph-con} of
Lemma~\ref{lem:finite-hypograph} are an exact reduction, rather than a
sampling or gridding approximation, and are the key finite step of the
game formulation.}
\rev{The direction of the curvature is essential: the follower
minimises a function that is concave in its own parameter, so an
extreme minimiser always exists. A convex margin may attain its minimum in the relative interior of $\mathcal{R}_h$, thus a vertex reduction is not generally valid in that case.}
To ease notation, we write the fixed attack vertices and the affine residuals as follows
\begin{equation*}
 w_s^{h,k}:=w_s^h(\rho^{h,k}),
 \  s=1,\ldots,q_{\Delta,h},
 \ k=1,\ldots,K_h,
\end{equation*}
\begin{align*}
 \Phi_{i,r,l,s}^{j,h,k}(b,u)
:=C_i\!\left(Av_r^j(b)+Bu+c+Ed_l+w_s^{h,k}\right)-b_i.
\end{align*}
For fixed $b\in\mathcal B_j$ and $(h,k)\in \mathcal{H}\times [1,K_h]$, we define explicitly the two
certificate feasible sets
\begin{align}
 \mathcal F_{\cI}^{j,h,k}(b):=
 \bigl\{(\eta,\{u_{r,l,s}\}):\;&
 \Phi_{i,r,l,s}^{j,h,k}(b,u_{r,l,s})\leq-\eta,
 \ \forall i,r,l,s,\nonumber\\[-1mm]
 &F_u u_{r,l,s}\leq g_u,\ \forall r,l,s\bigr\},
 \label{eq:F-ci}\\
 \mathcal F_{\rci}^{j,h,k}(b):=
 \bigl\{(\eta,\{u_r\}):\;&
 \Phi_{i,r,l,s}^{j,h,k}(b,u_r)\leq-\eta,
 \ \forall i,r,l,s,\nonumber\\[-1mm]
 &F_u u_r\leq g_u,\ \forall r\bigr\}.
 \label{eq:F-rci}
\end{align}

\begin{proposition}[Regional LP building block]
\label{prop:regional-LP}
Given any $(j,h)\in \mathcal{J}\times \mathcal{H}$, the controlled invariance regional game
$V_{\cI}^{j,h}$ is the optimal value of
\begin{subequations}\label{eq:ci-regional-LP}
\begin{align}
 \max_{b,t,\{u_{r,l,s,k}\}}\quad&t
\\[-1mm]
 \st\quad
 G_jb & \leq \beta_j,
 \\
 \Phi_{i,r,l,s}^{j,h,k}(b,u_{r,l,s,k})
 &\leq-t, \ \forall i,r,l,s,k,
 \\
 F_u u_{r,l,s,k}& \leq g_u,
 \ \forall r,l,s,k.
\end{align}
\end{subequations}
The robust game $V_{\rci}^{j,h}$ is the optimal value of
\begin{subequations}\label{eq:rci-regional-LP}
\begin{align}
 \max_{b,t,\{u_{r,k}\}}\quad&t
 \\[-1mm]
 \st\quad
 G_jb&\leq \beta_j,
 \\
 \Phi_{i,r,l,s}^{j,h,k}(b,u_{r,k})
 &\leq-t, \ \ \forall i,r,l,s,k,
 \\
 F_u u_{r,k}&\leq g_u,
 \ \ \forall r,k.
\end{align}
\end{subequations}
\end{proposition}

\begin{proof}
By Lemma~\ref{lem:finite-hypograph}, \revthree{it suffices to impose}
$t\leq\varepsilon_\nu^{j,h}(b,\rho^{h,k})$\revthree{; for fixed $k$,}
Proposition~\ref{prop:regional-certificates} and
\eqref{eq:F-ci}--\eqref{eq:F-rci} \revthree{give}
\begin{equation}
 \varepsilon_\nu^{j,h}(b,\rho^{h,k})
 =\max\{\eta:\exists\bm u^k,
 (\eta,\bm u^k)\in\mathcal F_\nu^{j,h,k}(b)\}.
 \label{eq:value-function-F}
\end{equation}
\revthree{The feasible set is downward closed in $\eta$: if $(\eta,\bm
u^k)$ is feasible, so is $(\tilde\eta,\bm u^k)$ for every
$\tilde\eta\leq\eta$. Since the maximum in \eqref{eq:value-function-F}
is attained,}
\begin{equation*}
 t\leq\varepsilon_\nu^{j,h}(b,\rho^{h,k})
 \Longleftrightarrow\quad
 \exists\bm u^k:\ (t,\bm u^k)\in\mathcal F_\nu^{j,h,k}(b).
 \end{equation*}
Stacking the $K_h$ equivalences yields
\eqref{eq:ci-regional-LP} and \eqref{eq:rci-regional-LP}.
\end{proof}

\revthree{Both programmes of Proposition~\ref{prop:regional-LP} are
LPs,} since $v_r^j(b)$ is affine in $b$ and $w_s^{h,k}$ is fixed. They provide the constraint blocks associated
with one pair $(\mathcal B_j,\mathcal R_h)$; the complete game on
$\mathcal B_j$ requires all $h\in\mathcal H$, as shown next.

\begin{remark}
\label{rem:rho-information}
In \eqref{eq:rci-regional-LP}, $u_{r,k}$ may depend on $\rho^{h,k}$ but
is shared across all $d_l$ and $w_s^{h,k}$, following the order
$\forall\rho\,\forall x\,\exists u\,\forall d\,\forall\delta$.
Replacing every $u_{r,k}$ by one $u_r$ imposes instead
$\forall x\,\exists u\,\forall\rho\,\forall d\,\forall\delta$,
for a controller derived before $\rho$, and possibly smaller
margin.
\end{remark}

Proposition~\ref{prop:regional-certificates} concerns a fixed pair
$(b,\rho)$, while Lemma~\ref{lem:finite-hypograph} and
Proposition~\ref{prop:regional-LP} fix $(j,h)$ and cover every $(b,\rho)\in\mathcal B_j\times\mathcal R_h$. In the game, however,
$S(b)$ must be safe against \emph{every} $\rho\in\mathcal R$, so we fix
$j$ and account for all $h\in\mathcal H$, setting
$V_{\nu,j}^\star:=\max_{b\in\mathcal B_j}\min_{\rho\in\mathcal R}
\varepsilon_\nu(b,\rho)$.

\begin{proposition}[One LP per defender region against all attacker regions]
\label{prop:one-LP-per-j}
For each $j\in\mathcal J$, $V_{\nu,j}^\star$ is obtained from one LP
containing the constraint blocks for all $h\in\mathcal H$ and
$k=1,\ldots,K_h$.
For controlled invariance, the LP is
\begin{subequations}\label{eq:ci-one-LP-j}
\begin{align}
 V_{\cI,j}^\star=\max_{b,t,\{u_{r,l,s,h,k}\}}\quad&t
 \\[-1mm]
 \st\quad
 G_jb&\leq \beta_j,
 \\
 \Phi_{i,r,l,s}^{j,h,k}(b,u_{r,l,s,h,k})
 &\leq-t, \ \  \forall h,k,i,r,l,s,
 \\
 F_u u_{r,l,s,h,k}&\leq g_u,
 \ \ \forall h,k,r,l,s.
\end{align}
\end{subequations}
For robust controlled invariance, it is
\begin{subequations}\label{eq:rci-one-LP-j}
\begin{align}
 V_{\rci,j}^\star=\max_{b,t,\{u_{r,h,k}\}}\quad&t
 \\ \st\quad
 G_jb&\leq \beta_j,
 \\
 \Phi_{i,r,l,s}^{j,h,k}(b,u_{r,h,k})
 &\leq-t, \ \ 
 \forall h,k,i,r,l,s,
 \\
 F_u u_{r,h,k}&\leq g_u,
\ \ \forall h,k,r.
\end{align}
\end{subequations}
\end{proposition}

\begin{proof}
Since $\mathcal R=\bigcup_{h\in\mathcal H}\mathcal R_h$, \rev{for
every $b\in\mathcal B_j$ we have
$\min_{\rho\in\mathcal R}\varepsilon_\nu(b,\rho)
=\min_{h\in\mathcal H}\min_{\rho\in\mathcal R_h}
\varepsilon_\nu^{j,h}(b,\rho)$.}
Lemma~\ref{lem:finite-hypograph} reduces the inner minimum on every
$\mathcal R_h$ to its $K_h$ vertices. A common lower bound $t$ must
therefore satisfy the certificate block associated with every pair
$(h,k)$. Stacking all blocks while retaining the same $b$ and
$t$ gives \eqref{eq:ci-one-LP-j} and \eqref{eq:rci-one-LP-j}.
\end{proof}

\begin{remark}[\rev{Size of the regional LPs}]
\label{rem:defender-region-LP-size}
\revtwo{Both LPs are feasible: any $b\in\mathcal B_j$ and admissible controls satisfy the successor inequalities for sufficiently negative $t$. Their optimal values are finite and attained because $\mathcal B_j$ and $U$ are compact and the constraints are closed.} \rev{\revthree{With} $\Sigma_1=\sum_hK_h$,
$\Sigma_2=\sum_hK_hq_{\Delta,h}$ and
$\Sigma_3=\sum_hK_h(pN_Dq_{\Delta,h}+m_U)$, the LP of
Proposition~\ref{prop:one-LP-per-j} on $\mathcal B_j$ has}
\begin{equation*}
\begin{aligned}
 N_{\cI,j}^{\rm var}&=p+1+n_uq_jN_D\Sigma_2,\\[-0.4mm]
 N_{\cI,j}^{\rm con}&=m_{B,j}+q_jN_D\Sigma_2(p+m_U),\\[-0.4mm]
 N_{\rci,j}^{\rm var}&=p+1+n_uq_j\Sigma_1,\quad
 N_{\rci,j}^{\rm con}=m_{B,j}+q_j\Sigma_3,
\end{aligned}
\end{equation*}
\rev{the regional block of Proposition~\ref{prop:regional-LP} being the
single term $h$. Each $(r,l,s,k)$ carries its own control under
controlled invariance, whereas under robust controlled invariance each
$(r,k)$ carries one control shared over $(l,s)$. The counts grow with all extreme parameter choices from all attacker
regions.}
\end{remark}

\begin{remark}
\label{rem:not-independent}
Let $g_{j,h}(b):=\min_{\rho\in\mathcal R_h}
\varepsilon_\nu^{j,h}(b,\rho)$. Solving
Proposition~\ref{prop:regional-LP} separately for each $h$ returns
$V_\nu^{j,h}=\max_{b\in\mathcal B_j}g_{j,h}(b)$, which permits a
different optimiser for every attacker region, whereas the order of
play requires one $b$ to face all of them,
\begin{equation*}
 V_{\nu,j}^\star=
 \max_{b\in\mathcal B_j}\min_{h\in\mathcal H}g_{j,h}(b)
 \ \leq\
 \min_{h\in\mathcal H}\max_{b\in\mathcal B_j}g_{j,h}(b),
\end{equation*}
the minimax inequality being possibly strict, since an optimiser for
one attacker region may have a poor margin on another. The values
$V_\nu^{j,h}$ are therefore insufficient to recover $V_{\nu,j}^\star$, and equality
would require additional structure that is not assumed. Proposition~\ref{prop:one-LP-per-j} preserves the
order by sharing $b$ and $t$ across all $h$ blocks, with only the scenario
controls being allowed to differ. 
\end{remark}
Thus, the LPs of distinct defender regions can be solved in
parallel, but within each LP the attacker blocks share $b$ and $t$ and
must be optimised jointly.
\rev{We next state the complete game and \revthree{its} solution
concept.}

\begin{definition}[\rev{Stackelberg set game and equilibrium}]
\label{def:stackelberg}
\rev{Let $\nu\in\{\cI,\rci\}$. In the \emph{Stackelberg set game} the
leader (defender) selects $b\in\mathcal B$, the follower (attacker)
observes $b$ and selects $\rho\in\mathcal R$, and the payoff
$\varepsilon_\nu(b,\rho)$ is maximised by the leader and minimised by
the follower. The follower best response correspondence and the leader guaranteed margin are}
\begin{equation}
 P_\nu(b):=\operatorname*{arg\,min}_{\rho\in\mathcal R}
 \varepsilon_\nu(b,\rho),
 \quad
 g_\nu(b):=\min_{\rho\in\mathcal R}\varepsilon_\nu(b,\rho).
 \label{eq:best-response}
\end{equation}
\rev{A pair $(b^\star,\rho^\star)\in\mathcal B\times\mathcal R$ is a
\emph{Stackelberg equilibrium} of the game if
$b^\star\in\arg\max_{b\in\mathcal B}g_\nu(b)$ and
$\rho^\star\in P_\nu(b^\star)$, and the associated \emph{Stackelberg value} is
$V_\nu^\star:=g_\nu(b^\star)=\varepsilon_\nu(b^\star,\rho^\star)$.}
\end{definition}

\rev{The game is zero sum, so every $\rho^\star\in P_\nu(b^\star)$
returns the same payoff and the Stackelberg value is unambiguous,}
\begin{equation}
 V_\nu^\star=
 \max_{b\in\mathcal B}\ \min_{\rho\in\mathcal R}
 \varepsilon_\nu(b,\rho),
 \qquad \nu\in\{\cI,\rci\}.
 \label{eq:global-game}
\end{equation}
\rev{Existence of an equilibrium is not assumed. It follows
constructively from Theorem~\ref{thm:global-finite-LPs} and
Proposition~\ref{prop:recovery} below, which also return $b^\star$ and
one $\rho^\star\in P_\nu(b^\star)$.}
The following theorem is the main result since it lifts the regional
elimination to the full defender and attacker parameter spaces without
changing the order of play.
\begin{theorem}[Exact solution of the set game]
\label{thm:global-finite-LPs}
Under Assumption~\ref{ass:polytopes} and the resulting validity-domain
decompositions \eqref{eq:partitions}--\eqref{eq:Dvertices}, the value of
\eqref{eq:global-game} is
\begin{equation}
 V_\nu^\star=
 \max_{j\in\mathcal J}V_{\nu,j}^\star,
 \qquad \nu\in\{\cI,\rci\},
 \label{eq:global-finite-result}
\end{equation}
where $V_{\cI,j}^\star$ and $V_{\rci,j}^\star$ are the LP values in
\eqref{eq:ci-one-LP-j} and \eqref{eq:rci-one-LP-j}, respectively.
\revthree{Hence} the complete Stackelberg set game is solved exactly by
$|\mathcal J|$ LPs\revthree{, each corresponding to one safe set vertex
region and containing all attack set vertex regions}.
\end{theorem}

\begin{proof}
By the finite cover
$\mathcal B=\bigcup_{j\in\mathcal J}\mathcal B_j$,
$V_\nu^\star=\max_{j\in\mathcal J}\max_{b\in\mathcal B_j}
\min_{\rho\in\mathcal R}\varepsilon_\nu(b,\rho)
=\max_{j\in\mathcal J}V_{\nu,j}^\star$.
Proposition~\ref{prop:one-LP-per-j} shows that each
$V_{\nu,j}^\star$ is the value of one LP. Taking the largest of this
finite collection establishes \eqref{eq:global-finite-result}.
\end{proof}
\revtwo{\revthree{Overlaps between validity domains do not affect the optimal
value in Theorem~\ref{thm:global-finite-LPs}}, since
\revthree{$\varepsilon_\nu(b,\rho)$} depends on the sets $S(b)$ and $\Delta(\rho)$, independently of their vertex representations. A parameter belonging to multiple domains therefore yields the same margin in each representation and does not alter the global maximisation. Moreover, \revthree{if} the safe set vertices are affine throughout
$\mathcal B$, the complete game reduces to one LP, \revthree{with}
additional attacker validity regions contributing only further
constraint blocks.}


\section{Optimal parameters, design, and computation}
\label{sec:subtleties}

\subsection{Recovery of optimal parameters and an active certificate tuple}
\revtwo{The \revthree{reformulation} recovers optimal defender and
attacker parameters \revthree{and} a state, control, disturbance and
attack tuple attaining an active successor constraint.}

\begin{proposition}[Recovery of optimal parameters and an active certificate tuple]
\label{prop:recovery}
Let $j^\star\in\arg\max_{j\in\mathcal J}V_{\nu,j}^\star$ and \revthree{$(b^\star,t^\star)$ optimise its LP in}
Proposition~\ref{prop:one-LP-per-j}. Then $t^\star=V_\nu^\star$.
With
\begin{equation}
 (h^\star,k^\star)\in
 \operatorname*{arg\,min}_{h\in\mathcal H,\,1\leq k\leq K_h}
 \varepsilon_\nu^{j^\star,h}(b^\star,\rho^{h,k}),
 \label{eq:recover-rho}
\end{equation}
$\rho^\star:=\rho^{h^\star,k^\star}$ \rev{satisfies
$\rho^\star\in P_\nu(b^\star)$, so $(b^\star,\rho^\star)$ is a
Stackelberg equilibrium per Definition~\ref{def:stackelberg}}, and
$\varepsilon_\nu(b^\star,\rho^\star)=t^\star$. Resolving the
certificate LP of
Proposition~\ref{prop:regional-certificates} at
$(b^\star,\rho^\star)$, there is an active tuple
$(i^\star,r^\star,l^\star,s^\star)$ with
\begin{equation}
 C_{i^\star}\!\left(Ax^\star+Bu^\star+c+Ed^\star+\delta^\star\right)
 =b_{i^\star}^\star-t^\star,
 \label{eq:critical-transition}
\end{equation}
where $x^\star=v_{r^\star}^{j^\star}(b^\star)$, $d^\star=d_{l^\star}$, $\delta^\star=w_{s^\star}^{h^\star}(\rho^\star)$ and $u^\star=u_{r^\star,l^\star,s^\star}$ for controlled invariance,
$u^\star=u_{r^\star}$ for the robust pattern. 
\revtwo{The tuple $(x^\star,u^\star,d^\star,\delta^\star)$ therefore activates a successor constraint at the optimal certificate margin.}
\end{proposition}

\begin{proof}
The first claim follows from Theorem~\ref{thm:global-finite-LPs}.
Lemma~\ref{lem:finite-hypograph}, applied on every attacker region,
shows that an attacker minimiser can be chosen among the finite pairs
$(h,k)$ in \eqref{eq:recover-rho}, and
$\varepsilon_\nu(b^\star,\rho^\star)=t^\star$ follows from optimality
of the common hypograph variable $t^\star$. At $(b^\star,\rho^\star)$, the optimal certificate margin is finite
and attained. At least one successor inequality must be active\rev{: if
every successor inequality were strictly slack, then $\varepsilon$
could be increased with the certificate controls unchanged, because
$\varepsilon$ does not enter the input constraints, contradicting
optimality}. Selecting any active inequality gives
\eqref{eq:critical-transition}.
\end{proof}

\begin{remark}[Nonuniqueness and numerical recovery]
\label{rem:nonuniqueness}
\revtwo{\revthree{T}he optimal parameters and active certificate tuple need not be unique. An active tuple need not minimise the margin after its associated control is reoptimised; consequently, activity alone does not identify a worst state and attack realisation.}
\end{remark}

\rev{Two orderings \revthree{help interpret a computed value}.}

\begin{corollary}[\rev{Monotonicity of the game value}]
\label{cor:monotonicity}
\rev{Under Assumption~\ref{ass:polytopes},
$V_{\rci}^\star\leq V_{\cI}^\star$. Moreover, if
$\mathcal R'\subseteq\mathcal R$ is a nonempty polytope and
$V_\nu^\star(\mathcal R')$ denotes the value of
\eqref{eq:global-game} with $\mathcal R$ replaced by $\mathcal R'$,
then $V_\nu^\star(\mathcal R')\geq V_\nu^\star$ for both $\nu\in\{\cI,\rci\}$.}
\end{corollary}

\begin{proof}
\rev{The pointwise ordering $\varepsilon_{\text{rci}}(b,\rho)\leq\varepsilon_{\text{ci}}(b,\rho)$ is preserved by the
minimisation over $\rho$ and then by the maximisation over $b$.
Shrinking the feasible set of the inner minimisation cannot decrease
$\min_{\rho}\varepsilon_\nu(b,\rho)$ for a fixed $b$, and the
subsequent maximisation preserves the inequality.}
\end{proof}

\begin{remark}[\rev{Stackelberg value and follower ties}]
\label{rem:stackelberg-value}
\rev{The set $P_\nu(b^\star)$ of Definition~\ref{def:stackelberg} need
not be a singleton, but the payoff is zero sum, so every selection from
it returns the same value: the optimistic and pessimistic Stackelberg
values coincide with \eqref{eq:global-game}
\cite{BasarOlsder1999}, no tie breaking rule is needed, and the nonuniqueness of $\rho^\star$ is immaterial. This is why the bilevel
problem \cite{SinhaMaloDeb2018} collapses in our case to the single level
programmes of Proposition~\ref{prop:one-LP-per-j}. The gap
$V_{\cI}^\star-V_{\rci}^\star\geq0$ is the value of observing
$(d,\delta)$ before committing to the control.}
\end{remark}

\subsection{Safe set design subject to a guaranteed margin}
\label{sec:design}

Maximising only $t$ may return a small set. Let $\bar\varepsilon\geq0$
be a prescribed margin and recall the hard state constraint set $\mrev{\mathcal X}$ from
\eqref{eq:hard-X} \rev{and the guaranteed margin $g_\nu$ from
\eqref{eq:best-response}}.

\begin{corollary}[Affine safe set size design]
\label{cor:size-design}
Suppose $J_{\rm size}(b)$ is affine on every $\mathcal B_j$ and let
$\kappa\geq0$. The problem
\begin{subequations}\label{eq:size-game}
\begin{align}
 \max_{b\in\mathcal B}\quad&
 J_{\rm size}(b)+\kappa g_\nu(b)
 \\
 \st\quad&g_\nu(b)\geq\bar\varepsilon,
 \qquad S(b)\subseteq \mrev{\mathcal X}
\end{align}
\end{subequations}
is solved by the same finite family of $|\mathcal J|$ LPs\rev{, in the
sense that its value is the largest of their optimal values, infeasible
regions being discarded}. By fixing a
defender validity region $\mathcal B_j$, we obtain its LP by replacing
the objective of \eqref{eq:ci-one-LP-j} or
\eqref{eq:rci-one-LP-j} with $J_{\rm size}(b)+\kappa t$ and adding
\begin{equation}
 t\geq\bar\varepsilon,\qquad
 H_Xv_r^j(b)\leq h_X,\quad r=1,\ldots,q_j.
 \label{eq:size-extra-constraints}
\end{equation}
\end{corollary}

\begin{proof}
By fixing $\mathcal B_j$, the first inequality in
\eqref{eq:size-extra-constraints} certifies
$g_\nu(b)\geq\bar\varepsilon$. The remaining inequalities are
equivalent to $S(b)\subseteq \mrev{\mathcal X}$ because they hold at every vertex of
$S(b)$. All added constraints and the regional objective are affine.
If $\kappa>0$, optimality forces $t=g_\nu(b)$. If $\kappa=0$, $t$
need not equal $g_\nu(b)$; it serves only as the required safety
certificate. Theorem~\ref{thm:global-finite-LPs} then applies without
modification.
\end{proof}

\rev{Unlike the margin game, \eqref{eq:size-game} may be infeasible on a
region, since $\bar\varepsilon$ and the hard constraint act jointly;
such a region is simply dropped from the comparison. If all
$|\mathcal J|$ LPs are infeasible, that is itself a certificate that
$(C,\mathcal B)$ cannot deliver the requested margin, and $\mathcal B$, the facet matrix $C$, or the admissible set $\mathcal
R$ must be revised.}
\revtwo{The components of $b$ specify the offsets of the defining
hyperplanes of $S(b)$. Thus, $J_{\mathrm{size}}(b)=\omega^\top b$,
with $\omega>0$, provides a linear surrogate for set size,
although increasing an offset associated with a redundant
inequality need not enlarge the set. Alternatively, under
the row normalisation and for $b\geq0$, maximising $\tau$
subject to $\tau\leq b_i$ for all $i$ maximises the radius of a norm ball centred at the origin and contained
in $S(b)$. Both formulations preserve the LP structure.}

\section{Numerical illustration: frequency regulation under additive attacks}
\label{sec:toy-example}

\begin{figure}[t]
 \centering
 \includegraphics[width=\columnwidth]{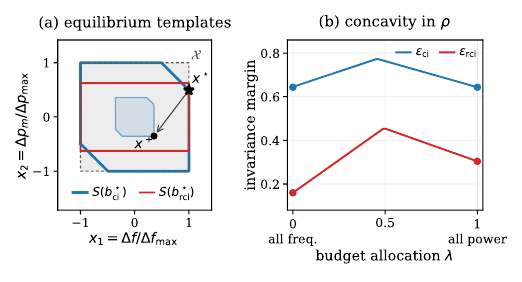}
 \caption{Stackelberg set game for the load frequency control
 model \eqref{eq:numerical-model} at $\Gamma=0.6$\revthree{, with
 $S(b)\subseteq\mathcal X$}. (a) Equilibrium
 templates for the two information patterns \revthree{within the
 envelope $\mathcal X$ (shaded, dashed)}, the contracted target
 $S_{V_{\cI}^\star}(b^\star_{\cI})$, and \revthree{an active certificate
 transition} of Proposition~\ref{prop:recovery} \revthree{(controlled
 invariance)}. (b) Margin along the full budget edge
 of $\mathcal R$, from the whole budget on the frequency channel
 ($\lambda=0$) to the whole budget on the power channel
 ($\lambda=1$); the map is concave and its minimum is attained at a
 vertex, as Lemma~\ref{lem:finite-hypograph} asserts.}
 \label{fig:lfc}
\end{figure}

\subsection{Model and set parameterisations}
We consider a simplified two-state model of frequency
regulation, motivated by standard load-frequency models
\cite{Bevrani2014}, with generator and battery actuation. The states represent frequency and mechanical power
deviations, normalised by $0.2$ Hz and $0.2$ p.u.,
respectively. The two control inputs are the generator
and battery commands, each normalised by $0.3$ p.u.,
while $d$ represents a load deviation normalised by  $0.05$ p.u. At a sampling period of $0.2$ s, we use the discrete-time model \eqref{eq:system} with $c=0$
and
{\small
\begin{equation}
\begin{aligned}
A&=\begin{bmatrix}
0.8975 & 0.7529\\
-0.1506 & 0.5361
\end{bmatrix},
&
B&=\begin{bmatrix}
0.3123 & 1.4416\\
0.5709 & -0.1249
\end{bmatrix},\\
E&=\begin{bmatrix}
-0.2403\\
0.0208
\end{bmatrix},
\end{aligned}
\label{eq:numerical-model}
\end{equation}
}
with $U=\{u:\|u\|_\infty\leq1\}$, thus $m_U=4$, and $D=\conv\{-1,1\}$,
thus $N_D=2$. The operational envelope is
$\mathcal X=\{x:\|x\|_\infty\leq1\}$. The open loop matrix has
eigenvalues $0.7168\pm0.2841\mathrm{j}$.
The defender template uses the six unit norm facet directions
$C=[\,e_1\ e_2\ -e_1\ -e_2\ \gamma\ -\gamma\,]^\top$ with
$\gamma=(1,1)/\sqrt2$, and the centrally symmetric offsets
$b=(a,c,a,c,\beta,\beta)$ with
\begin{equation}
 a,c\in[0.5,1.2],\qquad \beta\in[0.55,1.20].
 \label{eq:lfc-B}
\end{equation}
The attacker splits a budget $\Gamma$ between a frequency channel and a
power channel injection, the first being the cheaper to corrupt. With
$C_\Delta=[\,e_1\ -e_1\ e_2\ -e_2\,]^\top$ and
$\rho=(\alpha_1,\alpha_1,\alpha_2,\alpha_2)$,
$\mathcal R=\{\rho:\alpha_1,\alpha_2\geq0,\
\alpha_1+2\alpha_2\leq\Gamma\}$ and
$\Delta(\rho)=[-\alpha_1,\alpha_1]\times[-\alpha_2,\alpha_2]$.
The extreme full-budget allocations act exclusively on
the first or second component, yielding nonnested attack sets.

\subsection{Validity decomposition and regional LPs}

The diagonal facets are redundant exactly when
$\beta\geq(a+c)/\sqrt2$, so \eqref{eq:partitions} has $|\mathcal J|=2$ with
$\mathcal B_1=\mathcal B\cap\{\sqrt2\beta\geq a+c\}$,
$\mathcal B_2=\mathcal B\cap\{\sqrt2\beta\leq a+c\}$, and $S(b)$ a box with $q_1=4$ vertices on $\mathcal B_1$ and a corner
cut hexagon with $q_2=6$ vertices on $\mathcal B_2$, for instance
$v_1^2(b)=(a,\sqrt2\beta-a)$ and $v_2^2(b)=(\sqrt2\beta-c,c)$; both maps
are affine on their region, as in \eqref{eq:Svertices}. \revthree{The bounds \eqref{eq:lfc-B} are such that neither region
degenerates; on their common boundary, some generators coincide, while
the affine representations remain valid.}
\revthree{We incorporate $S(b)\subseteq\mathcal X$ into $\mathcal B$
by appending to $G_jb\leq\beta_j$ the vertex inequalities
$H_Xv_r^j(b)\leq h_X$, $r=1,\ldots,q_j$, of
\eqref{eq:size-extra-constraints}; they are affine on $\mathcal B_j$
and, under \eqref{eq:lfc-B}, amount to $a,c\leq1$.} The
attacker vertices are affine throughout $\mathcal R$, so
$|\mathcal H|=1$, $q_{\Delta,1}=4$ and
$\vtx(\mathcal R)=\{\rho^0,\rho^1,\rho^2\}
=\{0,(\Gamma,\Gamma,0,0),(0,0,\tfrac\Gamma2,\tfrac\Gamma2)\}$, \revthree{so $K_1=3$:} the
follower chooses from a continuum and
Lemma~\ref{lem:finite-hypograph} leaves three candidates. At $\Gamma=0.6$, after eliminating the template symmetry equalities,
the LP of Proposition~\ref{prop:one-LP-per-j} on $\mathcal B_2$ has
$292$ variables and \revthree{$1471$} inequalities for controlled invariance and
$40$ and \revthree{$967$} for the robust pattern, matching
Remark~\ref{rem:defender-region-LP-size} \revthree{with $m_{B,2}=31$}.

\subsection{Stackelberg equilibrium and an active certificate tuple}

At $\Gamma=0.6$, the regional LPs give
$V_{\mathrm{ci},1}^\star=0.4786$,
$V_{\mathrm{ci},2}^\star=0.6436$, and
$V_{\mathrm{rci},1}^\star=V_{\mathrm{rci},2}^\star=0.1597$.
Thus, the optimal controlled invariant template is
$b_{\mathrm{ci}}^\star=(1,1,1,1,1.0570,1.0570)$,
improving on the box region, whose margin is 26\% lower.
For robust controlled invariance, both regions attain
the same value; one optimiser is
$b_{\mathrm{rci}}^\star=(1,0.6247,1,0.6247,1.2,1.2)$.
Both sets lie in $\mathcal X$ and have positive margins,
certifying the respective invariance properties.
At $\rho^0,\rho^1,\rho^2$, the certificate margins are,
respectively, $(0.9299,0.6436,0.6436)$ for
$b_{\mathrm{ci}}^\star$ and $(0.6039,0.1597,0.3039)$
for $b_{\mathrm{rci}}^\star$.
Concavity and the strictly larger midpoint margin
$0.7633$ on $\operatorname{conv}\{\rho^1,\rho^2\}$
give $\mathcal P_{\mathrm{ci}}(b_{\mathrm{ci}}^\star)
=\{\rho^1,\rho^2\}$.
For the robust pattern, the unique minimising vertex
and concavity give
$\mathcal P_{\mathrm{rci}}(b_{\mathrm{rci}}^\star)=\{\rho^1\}$.
An active controlled-invariance certificate tuple at
$\rho^\star=\rho^1$ is
$x^\star=(1,0.4949)$, $u^\star=(-1,-1)$,
$d^\star=-1$, and $\delta^\star=(0.6,0)$.
It gives $x^+=(0.3564,-0.3521)$, with
$x_1^+=1-V_{\mathrm{ci}}^\star$, as shown in Fig.~1(a).
These active tuples need not be worst realisations
(Remark~4.2). The example demonstrates exact computation of optimal safe
sets and attack-set best responses through a finite family
of LPs: shape optimisation improves the margin when the
disturbance and attack are observed, whereas a box is already
optimal under the robust information pattern.

\color{black}

\end{document}